\documentclass[12pt,a4paper]{article}
\usepackage[T1]{fontenc}
\usepackage[latin1]{inputenc}
\usepackage[english]{babel}
\usepackage{amsmath}
\usepackage{amsfonts}
\usepackage{amssymb}
\usepackage{amsthm}
\usepackage{graphicx}
\usepackage[left=2cm,right=2cm,top=2cm,bottom=2cm]{geometry}

\begin{document}

\begin{center}
\textbf{\large On the Distribution of the Product and the Sum of Generalized Shifted Gamma Random Variables}
\par\end{center}{\large \par}

\begin{center}
Pushpa N. Rathie$^1$, Arjun K. Rathie$^2$ and Luan C. de S. M. Ozelim$^3$ 
\par\end{center}

\begin{center}
\textit{\emph{$^1$Department of Statistics, University
of Brasilia, Campus Universitário Darcy Ribeiro, Brasilia, DF, 70910-900, Brazil, email: pushpanrathie@yahoo.com}}\textit{ }
\par\end{center}
\begin{center}

\textit{\emph{$^2$Department of Mathematics, School of Mathematical and Physical Sciences, Central University of Kerala, Riverside Transit Campus, Padennakkad P.O. Nileshwar, Kasaragod - 671 328, Kerala, India, email: akrathie@rediffmail.com}}\textit{ }
\par\end{center}

\begin{center}
\textit{\emph{$^3$Department of Civil and Environmental
Engineering, University of Brasilia, Campus Universitário Darcy Ribeiro, Brasilia, DF, 70910-900, Brazil,
email: luanoz@gmail.com}}\textit{ }
\par\end{center}

\begin{abstract}
\noindent \emph{In general, while obtaining the probability density function of sums and products of shifted random variables, ordinary analytical methods such as Fourier and Mellin transforms tend to provide integrals which cannot be expressed in terms of ordinary Meijer G and H functions. This way, the need of defining new functions which easily enable one to write such integrals in a closed-form is inherent to the development of this area of statistical sciences. By generalizing the Mellin transform which defines the H function, a new function is established. A direct application of the so-called $\widehat{I}$ is discussed while developing the probability density function of the sum and the product of shifted generalized gamma random variables. Important special cases of the $\widehat{I}$ and their applications in science are also discussed in order to show the applicability of the function hereby defined. 
}

\medskip


\medskip


\end{abstract}

\section{Introduction}

The use of shifted distributions in the modelling of practical situations has grown considerably over the last years. Such growth is deeply related to the possibilities that arise by inserting a shifting - or translational - parameter into the definition of a given random variable. For example, the possibility of positioning the center of mass of a given distribution along its support enables a better description of the data analyzed. Also, while modelling measured data, shifting parameters make it possible to correct some experimental limitations such as measurements above/below a given threshold.

Even though shifted distributions have become more frequent in the scientific community, standard analytical functions and methods which would be capable of dealing with this kind of distributions fail in their intent. This way, a new approach to this problem is required.

In the present paper, a new function, hereby named $\widehat{I}$ function, which generalizes the well-known H function is defined.

In short, this new function changes the usual gamma functions in the contour integral representation by H functions themselves. This way, the $\widehat{I}$ function contemplates as special cases the H function \cite{Mathai}, the I function \cite{AKRathie}, the Y function \cite{YilmazY}, the Generalized Upper Incomplete Fox H function \cite{YilmazH} and other hypergeometric-type special functions.

In order to show the applicability of the $\widehat{I}$ function, the latter is used to represent the distributions of the sum and the product of generalized shifted gamma random variables. This result could not be achieved by means of the existent functions, which ensures the need of the new function defined in the present paper.

In the next section, the definition of the H function and its Mellin transform are presented which are used in the definition of the $\widehat{I}$ function.

\section{The H- function}

The H - function (see \cite{Springer},\cite{Mathai} and \cite{Braaksma}) is defined as an contour complex integral which
contain gamma functions in their integrands by 
\begin{gather}
H_{p,q}^{m,n}\left[z\;\bigg|\begin{array}{cccccc}
(a_{1},A_{1}), & \ldots, & (a_{n},A_{n}), & (a_{n+1},A_{n+1}), & \ldots, & (a_{p},A_{p})\\
(b_{1},B_{1}), & \ldots, & (b_{m},B_{m}), & (b_{m+1},B_{m+1}), & \ldots, & (b_{q},B_{q})
\end{array}\right]\nonumber \\
=\frac{1}{2\pi i}\int_{L}\frac{{\displaystyle \prod_{j=1}^{m}\Gamma(b_{j}+B_{j}s){\displaystyle \prod_{j=1}^{n}\Gamma(1-a_{j}-A_{j}s)}}}{{\displaystyle \prod_{j=m+1}^{q}\Gamma(1-b_{j}-B_{j}s){\displaystyle \prod_{j=n+1}^{p}\Gamma(a_{j}+A_{j}s)}}}z^{-s}ds,\label{eq:Hfunccontour}
\end{gather}
where \ $A_{j}$ and $B_{j}$ are assumed to be positive quantities
and all the $a_{j}$ and $b_{j}$ may be complex. The contour $L$ runs
from $c-i\infty$ to $c+i\infty$ such that the poles of $\Gamma(b_{j}+B_{j}s)$,
\ $j=1,\ldots,m$ lie to the left of $L$ and the poles of $\Gamma(1-a_{j}-A_{j}s)$,
$j=1,\ldots,n$ \ lie to the right of $L$.

The Mellin transform of the H -fuction is 
\begin{gather}
\int_{0}^{\infty}x^{s-1}H_{p,q}^{m,n}\left[cx\;\bigg|\begin{array}{cccccc}
(a_{p},A_{p})\\
(b_{q},B_{q})
\end{array}\right]dx=\frac{c^{-s}{\displaystyle \prod_{j=1}^{m}\Gamma(b_{j}+B_{j}s){\displaystyle \prod_{j=1}^{n}\Gamma(1-a_{j}-A_{j}s)}}}{{\displaystyle \prod_{j=m+1}^{q}\Gamma(1-b_{j}-B_{j}s){\displaystyle \prod_{j=n+1}^{p}\Gamma(a_{j}+A_{j}s)}}}.\label{eq:mellintransfH}
\end{gather}

Given these, one shall proceed to define the
generalized H function.

\section{The Generalized H function: The $\widehat{I}$ function}

The generalized H function, hereby named $\widehat{I}$ function, can be defined as a contour complex integral which
contain H functions in their integrands. For simplicity and in order to address a few problems of interest, the following simplified definition can be given 

\begin{gather}
\widehat{I}_{m}\left[z\;\bigg|\begin{array}{c}
(\overline{a}_{m,1},\hat{a}_{m,1},\overline{A}_{m,1}),(\overline{a}_{m,2},\hat{a}_{m,2},\overline{A}_{m,2}),(\overline{a}_{m,3},\hat{a}_{m,3},\overline{A}_{m,3})\\
(\overline{b}_{m,1},\hat{b}_{m,1},\overline{B}_{m,1}),(\overline{b}_{m,2},\hat{b}_{m,2},\overline{B}_{m,2})\\
\left(\overline{\gamma}_{m},\overline{\Gamma}_{m},\overline{\pi}_{m},\overline{\Pi}_{m},\overline{\rho}_{m},\overline{\sigma}_{m}\right)\\
\left(\overline{\alpha}_{m},\overline{\beta}_{m},\overline{\Lambda}_{m},\overline{\Theta}_{m},\overline{\zeta}_{m},\overline{\eta}_{m}\right)
\end{array}\right]\nonumber \\
=\frac{1}{2\pi i}\int_{L}\Upsilon\left(s\right)z^{-s}ds,\label{eq:GenHfunccontour}
\end{gather}
in which $\Upsilon(s)$ is the Mellin transform of the new function and can be explicitly given as:

\begin{equation}
\begin{array}{c}
\Upsilon(s)=\prod_{j=1}^{m}(\overline{\alpha}_{j}s+\overline{\beta}_{j})^{\overline{\Lambda}_{j}s+\overline{\Theta}_{j}}e^{\overline{\zeta}_{j}s+\overline{\eta}_{j}}\\
\left(H{}_{3,2}^{2,1}\left[\left(\overline{\gamma}_{j}+\overline{\Gamma}_{j}s\right)^{\overline{\pi}_{j}+\overline{\Pi}_{j}s}\;\bigg|\begin{array}{c}
(\overline{a}_{j,1}+\hat{a}_{j,1}s,\overline{A}_{j,1}),(\overline{a}_{j,2}+\hat{a}_{j,2}s,\overline{A}_{j,2}),(\overline{a}_{j,3}+\hat{a}_{j,3}s,\overline{A}_{j,3})\\
(\overline{b}_{j,1}+\hat{b}_{j,1}s,\overline{B}_{j,1}),(\overline{b}_{j,2}+\hat{b}_{j,2}s,\overline{B}_{j,2})
\end{array}\right]\right)^{\overline{\rho}_{j}s+\overline{\sigma}_{j}}
\end{array}\label{eq:MellinGenH}
\end{equation}

where $\overline{A}_{j,k}$ and $\overline{B}_{j,k}$ are assumed to
be positive real quantities, the $\overline{a}_{j,k}$, $\hat{a}_{j,k}$, $\overline{b}_{j,k}$, $\hat{b}_{j,k}$, $\overline{\Pi}_{j}$, $\overline{\pi}_{j}$,  $\overline{\Gamma}_{j}$, $\overline{\gamma}_{j}$,
$\overline{\rho}_{j}$, $\overline{\sigma}_{j}$, $\overline{\alpha}_{j}$, $\overline{\beta}_{j}$, $\overline{\Lambda}_{j}$, $\overline{\Theta}_{j}$, $\overline{\zeta}_{j}$, $\overline{\eta}_{j}$, $j=1,...,m$
are real numbers. The contour $L$ runs from $c-i\infty$ to $c+i\infty$, where $c$ is a real number,
and exists in accordance to Mellin inversion theorem, taking into account all the singularities.

Based on the definitions above, direct applications of the new function hereby introduced are shown in the next section. The new function can be further generalized by changing $H_{3,2}^{2,1}$ to $H_{p,q}^{m,n}$. This latter consideration is discussed in Section 6 of the present paper.

\section{Generalized Gamma Distribution and Distribution of its Product and Sum}

In the present section, a direct application of the generalized H function is developed. Let the probability density function of the generalized gamma distribution be given as:

\begin{equation}
f\left(x\right)=\frac{\gamma\beta^{\frac{\alpha}{\gamma}}}{\Gamma\left(\frac{\alpha}{\gamma}\right)}\left(x-\mu\right)^{\alpha-1}e^{-\beta\left(x-\mu\right)^{\gamma}},\; x>\mu;\;\alpha,\beta,\gamma>0\label{eq:pdfgenexp}
\end{equation}

One may notice that, in \eqref{eq:pdfgenexp}, by setting $\alpha=\gamma=k$; $\mu=0$ and $\beta=\lambda^{-k}$, the Weibull distribution with shape parameter $k$ and scale parameter $\lambda$ is retrieved. 
On the other hand, in \eqref{eq:pdfgenexp}, when $\alpha=\gamma=1$; $\mu=0$ and $\beta=\lambda$, the Exponential distribution with rate parameter $k$ is obtained. Also, by setting $\mu=0$, $\alpha=\gamma=2$ and $\beta=1/2\sigma^{2}$, the Rayleigh distribution with parameter $\sigma$ is recovered. By taking $\alpha=3$, $\mu=0$, $\gamma=2$ and $\beta=1/2a^{2}$ a Maxwell-Boltzmann distribution with parameter $a$ is obtained. Finally, by setting $\mu=0$, $\alpha=k$, $\gamma=1$ and $\beta=\lambda$, the generalized gamma distribution reduces to an Erlang distribution with shape parameter $k$ and rate parameter $\lambda$.
It is worth noticing that other gamma-type distributions may be expressed as special cases of \eqref{eq:pdfgenexp}.
The Generalized H function can be used to express \eqref{eq:pdfgenexp} as follows:

\newtheorem{theo}{Theorem}
\begin{theo}The probability density function of a generalized gamma random variable is given as:
\end{theo}

\begin{gather}
f(x)=\gamma\beta^{\frac{\alpha}{\gamma}}\mu^{\alpha-1} \widehat{I}_{1}\left[x\;\bigg|\begin{array}{c}
(1-\alpha,0,\gamma),\left(\frac{\alpha}{\gamma},0,0\right),(1,-1,0)\\
(0,0,1),(1-\alpha,-1,\gamma)\\
\left(\beta\mu^{\gamma},0,1,0,0,1\right)\\
\left(0,\mu,1,0,0,0\right)
\end{array}\right]\label{eq:InvMellinGenExpH}
\end{gather}

\begin{proof}
Consider the Mellin transform of the probability density function of a generalized gamma distribution. By means of \eqref{eq:pdfgenexp}, one may get:

\begin{equation}
\begin{array}{c}
M[f\left(x\right)](s)=\intop_{\mu}^{\infty}x^{s-1}\frac{\gamma\beta^{\frac{\alpha}{\gamma}}}{\Gamma\left(\frac{\alpha}{\gamma}\right)}\left(x-\mu\right)^{\alpha-1}e^{-\beta\left(x-\mu\right)^{\gamma}}dx\\
=\intop_{0}^{\infty}(y+\mu)^{s-1}\frac{\gamma\beta^{\frac{\alpha}{\gamma}}}{\Gamma\left(\frac{\alpha}{\gamma}\right)}y^{\alpha-1}e^{-\beta y^{\gamma}}dy
\end{array}\label{eq:Mellinpdfgenexp}
\end{equation}

On the other hand, consider the Mellin - Barnes representation
of the exponential function $e^{-\beta y^{\gamma}}$ : 
\begin{equation}
e^{-\beta y^{\gamma}}=\frac{1}{2\pi i}\int_{\mathcal{L}}\Gamma(s)[\beta y^{\gamma}]^{-s}ds,\label{eq:expcontour}
\end{equation}

This way, by inserting \eqref{eq:expcontour} into \eqref{eq:Mellinpdfgenexp},
the Mellin transform given in \eqref{eq:Mellinpdfgenexp} can be written
in terms of the H function as:

\begin{equation}
M[f(x)](s)=\gamma\beta^{\frac{\alpha}{\gamma}}\mu^{\alpha-1+s}H{}_{3,2}^{2,1}\left[\beta\mu^{\gamma}\;\bigg|\begin{array}{c}
(1-\alpha,\gamma),(\frac{\alpha}{\gamma},0),(1-s,0)\\
(0,1),(1-s-\alpha,\gamma)
\end{array}\right]\label{eq:MellinGenExpH}
\end{equation}

This way, by means of the inverse Mellin transform theorem, the representation \eqref{eq:InvMellinGenExpH} follows from \eqref{eq:MellinGenExpH}.
\end{proof}

\begin{theo}Consider the random variable $X=\underset{j=1}{\overset{N}{\prod}}X_{i}$ in which $X_{i}$, $i=1,..,N$ are independent generalized gamma random variables. This way, the probability density function of $X$ is given by:
\end{theo}

\begin{gather}
\begin{array}{c}
f_{X}(x)=\prod_{j=1}^{N}\left(\gamma_{j}\beta_{j}^{\frac{\alpha_{j}}{\gamma_{j}}}\mu_{j}^{\alpha_{j}-1}\right)\\
\widehat{I}_{N}\left[x\;\bigg|\begin{array}{c}
(1-\alpha_{1},0,\gamma_{1}),\left(\frac{\alpha_{1}}{\gamma_{1}},0,0\right),(1,-1,0)\\
(0,0,1),(1-\alpha_{1},-1,\gamma_{1})\\
\left(\beta_{1}\mu_{1}^{\gamma_{1}},0,1,0,0,1\right)\\
\left(0,\mu_{1},1,0,0,0\right)
\end{array};...;\begin{array}{c}
(1-\alpha_{N},0,\gamma_{N}),\left(\frac{\alpha_{N}}{\gamma_{N}},0,0\right),(1,-1,0)\\
(0,0,1),(1-\alpha_{N},-1,\gamma_{N})\\
\left(\beta_{N}\mu_{N}^{\gamma_{N}},0,1,0,0,1\right)\\
\left(0,\mu_{N},1,0,0,0\right)
\end{array}\right]
\end{array}\label{eq:ProductGenH}
\end{gather}

\begin{proof}
It is known that the Mellin transform of the distribution of the product
of independent random variables is the product of the Mellin transforms
of each variable \cite{Springer}, this way, the Mellin transform of the distribution
of the product of N independent generalized gamma random variables
is easily given by means of \eqref{eq:MellinGenExpH} as:

\begin{equation}
M[f_{X}(x)](s)=\prod_{j=1}^{N}\gamma_{j}\beta_{j}^{\frac{\alpha_{j}}{\gamma_{j}}}\mu_{j}^{\alpha_{j}-1+s}H{}_{3,2}^{2,1}\left[\beta_{j}\mu_{j}^{\gamma_{j}}\;\bigg|\begin{array}{c}
(1-\alpha_{j},\gamma_{j}),(\frac{\alpha_{j}}{\gamma_{j}},0),(1-s,0)\\
(0,1),(1-s-\alpha_{j},\gamma_{j})
\end{array}\right]\label{eq:MellinProducGenExpH}
\end{equation}

in which the subindex $j$ indicates the parameters of each distribution.
By means of \eqref{eq:GenHfunccontour} and \eqref{eq:MellinGenH}, the representation
in \eqref{eq:ProductGenH} easily follows.
\end{proof}

It is also worth noticing that the quotient of generalized gamma random variables can be easily obtained by making the substitution $s=2-s$ in \eqref{eq:MellinGenExpH} for the random variable which is in the denominator of the ratio \cite{Springer}. Due to the simplicity of the procedure, its full development is not presented in the present paper.

\begin{theo}Consider the random variable $Y=\underset{j=1}{\overset{N}{\sum}}X_{i}$ in which $X_{i}$, $i=1,..,N$ are independent generalized gamma random variables. This way, the probability density function of $Y$ is given by:
\end{theo}

\begin{gather}
\begin{array}{c}
f_{Y}(x)=\prod_{j=1}^{N}\left(\gamma_{j}\beta_{j}^{\frac{\alpha_{j}}{\gamma_{j}}}\right)\\
\widehat{I}_{N}\left[e^{x}\;\bigg|\begin{array}{c}
(1-\alpha_{1},0,\gamma_{1}),\left(\frac{\alpha_{1}}{\gamma_{1}},0,0\right),(1,0,0)\\
(0,0,1),(1,0,0)\\
\left(0,\beta_{1}^{1/\gamma_{1}},\gamma_{1},0,0,1\right)\\
\left(1,0,0,-\alpha_{1},-\mu_{1},0\right)
\end{array};...;\begin{array}{c}
(1-\alpha_{N},0,\gamma_{N}),\left(\frac{\alpha_{N}}{\gamma_{N}},0,0\right),(1,0,0)\\
(0,0,1),(1,0,0)\\
\left(0,\beta_{N}^{1/\gamma_{N}},\gamma_{N},0,0,1\right)\\
\left(1,0,0,-\alpha_{N},-\mu_{N},0\right)
\end{array}\right]
\end{array}\label{eq:SumGenH}
\end{gather}

\begin{proof}
In order to obtain the probability density function of the sum of
N independent shifted generalized gamma random variables, consider the
Laplace transform of the probability density function given in \eqref{eq:pdfgenexp}

\begin{equation}
\begin{array}{c}
L[f\left(x\right)](s)=\intop_{\mu}^{\infty}e^{-sx}\frac{\gamma\beta^{\frac{\alpha}{\gamma}}}{\Gamma\left(\frac{\alpha}{\gamma}\right)}\left(x-\mu\right)^{\alpha-1}e^{-\beta\left(x-\mu\right)^{\gamma}}dx\\
=\intop_{0}^{\infty}e^{-s(y+\mu)}\frac{\gamma\beta^{\frac{\alpha}{\gamma}}}{\Gamma\left(\frac{\alpha}{\gamma}\right)}y^{\alpha-1}e^{-\beta y^{\gamma}}dy
\end{array}\label{eq:Laplacepdfgenexp}
\end{equation}

On the other hand, consider the Mellin - Barnes representation
of the exponential function $e^{-\beta y^{\gamma}}$ , given in \eqref{eq:expcontour}. 

This way, by inserting \eqref{eq:expcontour} into \eqref{eq:Laplacepdfgenexp},
the Laplace transform given in \eqref{eq:Laplacepdfgenexp} can be written
in terms of the H function as:

\begin{equation}
L[f(x)](s)=\gamma\beta^{\frac{\alpha}{\gamma}}e^{-s\mu}s^{-\alpha}H{}_{3,2}^{2,1}\left[\beta s^{-\gamma}\;\bigg|\begin{array}{c}
(1-\alpha,\gamma),(\frac{\alpha}{\gamma},0),(1,0)\\
(0,1),(1,0)
\end{array}\right]\label{eq:LaplaceGenExpH}
\end{equation}

It is known that the Laplace transform of the distribution of the
sum of independent random variables is the product of the Laplace
transforms of each random variable \cite{Springer}, this way, the Laplace transform of the
distribution of the sum of N independent generalized gamma random
variables is easily given as:

\begin{equation}
L[f_{Y}(x)](s)=\prod_{j=1}^{N}\gamma_{j}\beta_{j}^{\frac{\alpha_{j}}{\gamma_{j}}}e^{-s\mu_{j}}s^{-\alpha_{j}}H{}_{3,2}^{2,1}\left[\beta_{j}s^{-\gamma_{j}}\;\bigg|\begin{array}{c}
(1-\alpha_{j},\gamma_{j}),(\frac{\alpha_{j}}{\gamma_{j}},0),(1,0)\\
(0,1),(1,0)
\end{array}\right]\label{eq:LaplaceSumGenExpH}
\end{equation}

By means of \eqref{eq:GenHfunccontour} and \eqref{eq:MellinGenH}, the representation
in \eqref{eq:SumGenH} easily follows.
\end{proof}

The formulas obtained for the sum of N independent generalized gamma random variables easily enable one to obtain the distribution of the linear combination of random variables of this type. This can be achieved by noticing that, if $X$ is a generalized gamma random variable with parameters $\alpha$, $\beta$, $\gamma$ and $\mu$, for a given constant value $A$, the random variable  $AX$ is also a generalized gamma random variable with parameters $\alpha$, $\beta/A^{\gamma}$, $\gamma$ and $A\mu$, respectively.

One may also notice that the results presented in the present section are a generalization of the ones in \cite{YilmazH} and \cite{YilmazY}. Also the new results are entirely dependent on the definition of the $\widehat{I}$ function.

\section{Special Cases of the $\widehat{I}$ function}

In the present section, a few special cases of the generalized H function are discussed.

\subsection{Standard H function}
The key equation to relate the new function hereby defined to the standard H function is the contour integral representation of the gamma function, given as:

\begin{equation}
\Gamma\left(z\right)=\frac{1}{2\pi i}\int_{\mathcal{L}}\frac{\Gamma(s)\Gamma(z)}{\Gamma(1+s)}ds,\label{eq:StdGamma}
\end{equation}

Thus, from \eqref{eq:Hfunccontour}, one shall get:

\begin{equation}
\Gamma\left(z\right)=H{}_{3,2}^{2,1}\left[1\;\bigg|\begin{array}{c}
(0,0),(1,1),(1,0)\\
(0,1),(z,0)
\end{array}\right],\label{eq:StdrGammaHFunc}
\end{equation}

On the other hand, by setting $\overline{\alpha}_{j}=\overline{\zeta}_{j}=\overline{\eta}_{j}=\overline{\Gamma}_{j}=\overline{\Pi}_{j}=\overline{\Lambda}_{j}=\overline{\rho}_{j}=0,\overline{\gamma}_{j}=\overline{\Theta}_{j}=\overline{\beta}_{j}=\overline{\pi}_{j}=\overline{\sigma}_{j}=1$,
and taking into account \eqref{eq:StdrGammaHFunc} the first produtory
present in the numerator of \eqref{eq:mellintransfY} is recovered.
The other produtories are easily obtained by similar procedures, ultimately
showing that the H function is a
special case of the $\widehat{I}$ function.

\subsection{I - function}
In \cite{AKRathie}, a generalized hypergeometric function has been defined. The main advantage of the so-called I function is that powers of gamma functions are considered inside the Mellin transform used to define such a function. Explicitly, in \cite{AKRathie} it has been defined the I function as the inverse Mellin transform of
the following function $I(s)$:

\begin{gather}
I\left(s\right)=\frac{{\displaystyle \prod_{j=1}^{m}\Gamma^{B_{j}}(b_{j}+\beta_{j}s){\displaystyle \prod_{j=1}^{n}\Gamma^{A_{j}}(1-a_{j}-\alpha_{j}s)}}}{{\displaystyle \prod_{j=m+1}^{q}\Gamma^{B_{j}}(1-b_{j}-\beta_{j}s){\displaystyle \prod_{j=n+1}^{p}\Gamma^{A_{j}}(a_{j}+\alpha_{j}s)}}}\label{eq:mellintransfY-1-1}
\end{gather}

This way, by setting $\overline{\alpha}_{j}=\overline{\zeta}_{j}=\overline{\eta}_{j}=\overline{\Gamma}_{j}=\overline{\Pi}_{j}=\overline{\Lambda}_{j}=\overline{\rho}_{j}=0,\overline{\gamma}_{j}=\overline{\Theta}_{j}=\overline{\beta}_{j}=\overline{\pi}_{j}=1,\overline{\sigma}_{j}=B_{j}$,
and taking into account \eqref{eq:StdrGammaHFunc} the first produtory
present in the numerator of \eqref{eq:mellintransfY-1-1} is recovered.
The other produtories are easily obtained by similar procedures, ultimately
showing that the I function is a special case of the $\widehat{I}$ function.

It is interesting to notice that recently the I function has found important application in the study of wireless communication systems.
In \cite{Ansari1} it has been derived the probability density function and cumulative distribution function of the sum of L independent but not necessarily identically distributed Gamma variates. Their result is applicable to the modelling of the output statistics of maximal ratio combining (MRC) receiver operating over Nakagami-m fading channels. The results are derived in terms of Meijer G function, H function and I function.
Following a similar approach as in \cite{Ansari1}, in \cite{Ansari2} it has been derived the probability density function and cumulative distribution function of the sum of L independent but not
necessarily identically distributed squared $\eta - \mu$ variates. As in the case of \cite{Ansari1}, the results are applicable to the output statistics of maximal ratio combining (MRC) receiver operating over $\eta - \mu$ fading channels. The I function has a major role in the development of the results of both \cite{Ansari1} and \cite{Ansari2}.

\subsection{Y - function}
In a recent paper \cite{YilmazY}, it has been proposed a generalization of the H function. In short, the so-called Y function generalizes the Mellin transform of the H function by changing the ordinary gamma function in \eqref{eq:mellintransfH} by Tricomi hypergeometric functions, defined as the following contour integral:

\begin{equation}
U\left(a,b,z\right)=\frac{1}{\Gamma\left(a\right)\Gamma\left(a-b+1\right)}\frac{1}{2\pi i}\int_{\mathcal{L}}\Gamma(s)\Gamma(1-b+s)\Gamma(a-s)z^{-s}ds,\label{eq:TricomiHyper}
\end{equation}

This way, by pre-multiplying the Tricomi hypergeomtric functions by
a constant raised to a linear function of $s$, in \cite{YilmazY} it has been defined the Y function as the inverse Mellin transform of the following function $Y(s)$:

\begin{gather}
\begin{array}{c}
Y\left(s\right)=\\
\frac{{\displaystyle \prod_{j=1}^{m}B_{j}^{\phi_{j}+b_{j}+\beta_{j}s-1}U(\phi_{j},\phi_{j}+b_{j}+\beta_{j}s,B_{j}){\displaystyle \prod_{j=1}^{n}A_{j}^{\varphi_{j}-a_{j}-\alpha_{j}s}U(\varphi_{j},\varphi_{j}-a_{j}-\alpha_{j}s+1,A_{j})}}}{{\displaystyle \prod_{j=n+1}^{p}A_{j}^{\varphi_{j}+a_{j}+\alpha_{j}s-1}U(\varphi_{j},\varphi_{j}+a_{j}+\alpha_{j}s,A_{j})}\underset{j=m+1}{\overset{q}{\prod}}B_{j}^{\phi_{j}-b_{j}-\beta_{j}s}U(\phi_{j},\phi_{j}-b_{j}-\beta_{j}s+1,B_{j})}
\end{array}\label{eq:mellintransfY}
\end{gather}

One may notice that \eqref{eq:mellintransfY} is misprinted in \cite{YilmazY}.

In order to present the relation between both functions, by means
of \eqref{eq:Hfunccontour}, it is easy to notice that:

\begin{equation}
U\left(a,b,z\right)=H{}_{3,2}^{2,1}\left[z\;\bigg|\begin{array}{c}
(1-a,1),(a,0),(a-b+1,0)\\
(0,1),(1-b,1)
\end{array}\right],\label{eq:TricomiHyperHfunc}
\end{equation}

Finally, by setting $\overline{\alpha}_{j}=\overline{\zeta}_{j}=\overline{\eta}_{j}=\overline{\Gamma}_{j}=\overline{\Pi}_{j}=\overline{\rho}_{j}=0,\overline{\beta}_{j}=\overline{\gamma}_{j}=B_{j},\overline{\Lambda}_{j}=\beta_{j},\overline{\Theta}_{j}=\phi_{j}+b_{j}-1,\overline{\pi}_{j}=\overline{\sigma}_{j}=1$
and by means of \ref{eq:TricomiHyperHfunc}, the first produtory present
in the numerator of \eqref{eq:mellintransfY} is recovered. The other
produtories are easily obtained by similar procedures, ultimately
showing that the Y - function is a special case of the $\widehat{I}$ function.

The Y function has found interesting applications in wireless communication systems. In special, \cite{YilmazY} obtained the PDF and CDF of a shifted gamma random variable in terms of the Y function. By means of their results, it was possible to simulate the outage capacity of a multicarrier transmission system through a slow Nakagami-m fading channel.

\subsection{Generalized Upper Incomplete Fox H function}
Also in recent paper \cite{YilmazH}, it has been proposed another function which generalize the H function. The so-called Generalized Upper Incomplete Fox H function generalizes the Mellin transform of the H function by changing the ordinary gamma functions in \eqref{eq:mellintransfH} by upper incomplete gamma functions, which may be defined by means of the following contour integral:

\begin{equation}
\Gamma\left(a,z\right)=\frac{1}{2\pi i}\int_{\mathcal{L}}\frac{\Gamma(s)\Gamma(s+a)}{\Gamma(1+s)}z^{-s}ds\label{eq:UpperGamma}
\end{equation}

In \cite{YilmazH} it has been defined the Generalized Upper Incomplete H function as the inverse Mellin transform of the following function $UI(s)$:

\begin{gather}
UI\left(s\right)=\frac{{\displaystyle \prod_{j=1}^{m}\Gamma(b_{j}+\beta_{j}s,B_{j}){\displaystyle \prod_{j=1}^{n}\Gamma(1-a_{j}-\alpha_{j}s,A_{j})}}}{{\displaystyle \prod_{j=m+1}^{q}\Gamma(1-b_{j}-\beta_{j}s,B_{j}){\displaystyle \prod_{j=n+1}^{p}\Gamma(a_{j}+\alpha_{j}s,A_{j})}}}\label{eq:mellintransfY-1}
\end{gather}

In order to present the relation between both functions, by means
of \eqref{eq:Hfunccontour}, it is easy to notice that:

\begin{equation}
\Gamma\left(a,z\right)=H{}_{3,2}^{2,1}\left[z\;\bigg|\begin{array}{c}
(0,0),(1,1),(1,0)\\
(0,1),(a,1)
\end{array}\right],\label{eq:UpperGammaHFunc}
\end{equation}

Thus, by setting $\overline{\alpha}_{j}=\overline{\zeta}_{j}=\overline{\eta}_{j}=\overline{\Gamma}_{j}=\overline{\Pi}_{j}=\overline{\Lambda}_{j}=\overline{\rho}_{j}=0,\overline{\gamma}_{j}=B_{j},\overline{\Theta}_{j}=\overline{\beta}_{j}=\overline{\pi}_{j}=\overline{\sigma}_{j}=1$,
and taking into account \eqref{eq:UpperGammaHFunc} the first produtory
present in the numerator of \eqref{eq:mellintransfY-1} is recovered.
The other produtories are easily obtained by similar procedures, ultimately
showing that the Generalized Upper Incomplete Fox H function is a
special case of the $\widehat{I}$ function.

As in the cases of both I and Y function, the Generalized Upper Incomplete H function has found useful applications in the modelling of wireless communication systems. In \cite{YilmazH} the PDF and CDF of the product of shifted exponential random variables has been obtained in terms of the Generalized Upper Incomplete H function. By applying such results, it has been possible to model the outage capacity of a multicarrier transmission system through a slow Rayleigh fading channel.

\section{Further generalization of $\widehat{I}$ function}

A generalization of the $\widehat{I}$ function, the $\hat{\hat{I}}$ function, can be defined as follows:

\begin{gather}
\hat{\hat{I}}{}_{m}\left[z\;\bigg|\begin{array}{c}
\{(\overline{a}_{m,1},\hat{a}_{m,1},\overline{A}_{m,1}),...,(\overline{a}_{m,\overline{p}_{m}},\hat{a}_{m,\overline{p}_{m}},\overline{A}_{m,\overline{p}_{m}}),(\overline{m}_{m},\overline{n}_{m},\overline{p}_{m},\overline{q}_{m})\}\\
\{(\overline{b}_{m,1},\hat{b}_{m,1},\overline{B}_{m,1}),...,(\overline{b}_{m,\overline{q}_{m}},\hat{b}_{m,\overline{q}_{m}},\overline{B}_{m,\overline{q}_{m}})\}\\
\{\left(\overline{\gamma}_{m},\overline{\Gamma}_{m},\overline{\pi}_{m},\overline{\Pi}_{m},\overline{\rho}_{m},\overline{\sigma}_{m}\right)\}\\
\{\left(\overline{\alpha}_{m},\overline{\beta}_{m},\overline{\Lambda}_{m},\overline{\Theta}_{m},\overline{\zeta}_{m},\overline{\eta}_{m}\right)\}
\end{array}\right]\nonumber \\
=\frac{1}{2\pi i}\int_{L}\overline{\Upsilon}\left(s\right)z^{-s}ds,\label{eq:GenHfunccontourTot}
\end{gather}
in which $\overline{\Upsilon}(s)$ is the Mellin transform of the new function and can be explicitly given as:

\begin{equation}
\begin{array}{c}
\overline{\Upsilon}(s)=\prod_{j=1}^{m}(\overline{\alpha}_{j}s+\overline{\beta}_{j})^{\overline{\Lambda}_{j}s+\overline{\Theta}_{j}}e^{\overline{\zeta}_{j}s+\overline{\eta}_{j}}\\
\left(H{}_{\overline{p}_{j},\overline{q}_{j}}^{\overline{m}_{j},\overline{n}_{j}}\left[\left(\overline{\gamma}_{j}+\overline{\Gamma}_{j}s\right)^{\overline{\pi}_{j}+\overline{\Pi}_{j}s}\;\bigg|\begin{array}{c}
(\overline{a}_{j,1}+\hat{a}_{j,1}s,\overline{A}_{j,1}),...,(\overline{a}_{j,\overline{p}_{j}}+\hat{a}_{j,\overline{p}_{j}},\overline{A}_{j,\overline{p}_{j}})\\
(\overline{b}_{j,1}+\hat{b}_{j,1}s,\overline{B}_{j,1}),...,(\overline{b}_{j,\overline{q}_{j}}+\hat{b}_{j,\overline{q}_{j}}s,\overline{B}_{j,\overline{q}_{j}})
\end{array}\right]\right)^{\overline{\rho}_{j}s+\overline{\sigma}_{j}}
\end{array}\label{eq:MellinGenHTot}
\end{equation}

where, as discussed in Section 3, $\overline{A}_{j,k}$ and $\overline{B}_{j,k}$ are assumed to
be positive real quantities, the $\overline{a}_{j,k}$, $\hat{a}_{j,k}$, $\overline{b}_{j,k}$, $\hat{b}_{j,k}$, $\overline{\Pi}_{j}$, $\overline{\pi}_{j}$,  $\overline{\Gamma}_{j}$, $\overline{\gamma}_{j}$,
$\overline{\rho}_{j}$, $\overline{\sigma}_{j}$, $\overline{\alpha}_{j}$, $\overline{\beta}_{j}$, $\overline{\Lambda}_{j}$, $\overline{\Theta}_{j}$, $\overline{\zeta}_{j}$, $\overline{\eta}_{j}$, $j=1,...,m$
are real numbers. The contour $L$ runs from $c-i\infty$ to $c+i\infty$, where $c$ is a real number,
and exists in accordance to Mellin inversion theorem, taking into account all the singularities.
The parameters $\overline{m}_{j},\overline{n}_{j},\overline{p}_{j}$ and $\overline{q}_{j}$, $j=1,...,m$, are positive integers following H function's definition.
It is important to notice that the braces sign used in the definition of the $\hat{\hat{I}}$ function denotes a vector whose components are to be written as the generic case shown.

\section{Conclusions}

In spite of the growth in the use of shifted random variables, analytical functions and methods are unable to give closed-form representations for sums and products of this type of random variables. 

In the present paper, in order to approach this issue, a new $\widehat{I}$ function, is defined. This new function generalizes the H function by changing the gamma functions present in its contour integral definition by H functions themselves. This way, not only the H function but also the I, Y and Generalized Upper Incomplete Fox H functions are shown to be special cases of the $\widehat{I}$ function.

The new function has been successfully applied to obtain the distributions of the sum and the product of generalized shifted gamma random variables, showing its value both in theoretical and practical backgrounds.

\bibliographystyle{amsplain}
\bibliography{RathieRathieOzelimIhat2013Final}

\end{document}